\documentclass{amsart}
\usepackage{times}

\issueinfo{00}
  {}
  {}
  {1995}

\newtheorem*{AM}{Axiom of Monotonicity}
\newtheorem*{AF}{Axiom of Fusion}
\newtheorem*{CT}{Continuum Theorem}
\newtheorem*{GCH}{Generalized Continuum Hypothesis (GCH)}
\newtheorem*{AC}{Axiom of Choice (AC)}
\newtheorem*{CBT}{Combinatorial Theorem}
\newtheorem*{UT}{Unification Theorem}

\theoremstyle{definition}

\theoremstyle{remark}

\numberwithin{equation}{section}

\begin{document}

\title{The Essence of Intuitive Set Theory}

\author{K. K. Nambiar}
\address{Formerly, Jawaharlal Nehru University, New Delhi, 110067, India}
\curraddr{1812 Rockybranch Pass, Marietta, Georgia, 30066-8015}
\email{nambiar@mediaone.net}


\date{June 25, 2001}



\begin{abstract}
Intuitive set theory is defined as the theory we get when we add 
the axioms, Monotonicity and Fusion, to ZF theory. Axiom of 
Monotonicity makes the Continuum Hypothesis true, and the Axiom 
of Fusion splits the unit interval into infinitesimals.
\end{abstract}

\maketitle

\section{INTRODUCTION}

\noindent The primary purpose of this paper is to give a clear 
definition of intuitive set theory (IST), 
so that researchers have all the 
necessary background to investigate the consistency 
of the two axioms that define
IST \cite{BePu:PM,Cantor:CFTTN,Nam:VIST}. G\"odel 
tells us that even though we 
are not in a position to  prove the consistency of a significant theory, 
we can prove its inconsistency, if it is inconsistent.
The secondary purpose of this paper is to explain IST
to the novice who has a passing acquaintance with the transfinite
cardinals of Cantor. 

\section{SEQUENCES AND SETS}

We will accept the fact that every number in the open interval $(0,1)$ can be 
represented \textit{uniquely} by an 
\textit{infinite} nonterminating binary sequence. 

For example, the infinite binary sequence 
\[
.10101010 \cdots
\]
can be recognized as the representation for the 
number $2/3$ and 
\[
.10111111 \cdots
\]
for the number $3/4$.
This in turn implies that an infinite subset of positive 
integers can be used to represent the numbers in the interval $(0,1)$.
Thus we have the set
\[
\{1,3,5,7,\cdots\}^{+}
\]
also as a representation for $2/3$. A binary
sequence that goes towards the right as 
above, we will call a \textit{right-sequence} 
and the corresponding set a \textit{right-set}, to make provision for  
a \textit{left-sequence} and a \textit{left-set}.
It is easy to see that the left sequence 
\[
\cdots 000010011.
\]
and the corresponding left-set 
\[
{}^{-}\{4,1,0\}
\]
can be used to represent the number $19$. In general, any nonnegative 
integer can be represented by a left-sequence, which eventually ends 
up in $0$s. Two's complement number system allows us to use 
left-sequences which eventually end up in $1$s to represent negative 
integers. Thus, we have the left-sequence, 
\[
\cdots 111100101.
\]
and the left-set
\[
{}^{-}\{\cdots 8,7,6,5,2,0\}
\]
representing the negative number $-27$. Adding up all these facts, we can 
claim that a two-way sequence can be used to represent any number 
on a real line. For example, the sequence, 
\[
\cdots 000010011.10101010 \cdots 
\]
and the corresponding two-way set
\[
{}^{-}\{4,1,0:1,3,5,7\cdots\}^{+}
\]
represent the number $19.6666 \cdots$. Similarly, the complement 
of this sequence, 
\[
\cdots 111101100.01010101 \cdots 
\]
and the corresponding two-way set 
\[
{}^{-}\{\cdots 8,7,6,5,3,2:2,4,6,8 \cdots\}^{+}
\]
represent the negative number $-19.6666 \cdots$. Note the restriction 
in our definition of a real number: the left sequence must eventually 
end up in either $1$s or $0$s. The number system we get when we put
no restriction on both the left-sequence and the right-sequence,
we will call the  
\textit{universal number system} (UNS). A universal number, whose 
left sequence is \textit{not} eventually-periodic, we will call a 
\textit{supernatural number}. The connection between 
the transcendental and supernatural numbers is explained next.

\section{UNIVERSAL NUMBER SYSTEM}

We will first explain why we have excepted eventually-periodic 
left-sequences from our definition of supernatural numbers. Consider 
the left-sequence 
\[
\cdots 101101101001001.
\]
with a periodic part $a=101=5$ of length $l_p=3$
and a nonperiodic part $b=001001=9$ of length $l_n=6$.
We can write the sequence formally as 
\[
b+\frac{a2^{l_n}}{1-2^{l_p}}
\]
which when evaluated gives
\[
-\frac{257}{7}.
\]
From this we infer that eventually-periodic left-sequences corresponds to 
negative rational numbers. A similar argument shows that 
eventually-periodic right-sequences represent positive rational numbers.

We want to show that corresponding to every transcendental number 
there is a supernatural number. 
Given a universal number $a$, the number we get when we flip the two-way 
infinite string around the binary point, we will write as $a^{F}$.
Consider the transcendental number 
\[\frac{\pi}{4}= \cdots000.11001000110\cdots \]
and
\[\Bigl({\frac{\pi}{4}}\Bigr)^{F}= \cdots01100010011.000\cdots \]
which gives the appearance of a number above all natural 
numbers. It is for this reason, we have 
called it a supernatural number, but of course, it is no
more supernatural than the transcendental number is transcendental.
From this example, we infer that corresponding to every
transcendental number in the interval $(0,1)$, there is a supernatural 
number. More generally, we can say that every irrational number 
in the interval $(0,1)$ has a corresponding supernatural number. 
By definition, an \textit{infinite recursive} subset of positive integers, 
is an infinite right-set with
a clear algorithm for its generation. The 
corresponding number in the interval $(0,1)$ is called a 
computable number. It is known from recursive function theory that
the cardinality of the set $R$ of these computable numbers is $\aleph_0$.
A number in the interval $(0,1)$, which is 
not computable, we will call an \textit{illusive} number. We will have 
more to say about irrational computable numbers, but before that
we want to take a cursory look at the transfinite cardinals of 
Cantor.

\section{TRANSFINITE CARDINALS}

Recall that every natural number can be represented by a set 
as given below. 

\begin{align*}
\{\}&=0, \\
\{0\}&=1, \\
\{0,1\}&=2, \\
\{0,1,2\}&=3, \\
&\vdots
\end{align*}
The advantage with this method is that we get an elegant way of
defining the first transfinite cardinal of Cantor, as
\[
{\aleph_0}=\{0,1,2,3,\cdots \}.
\]
The set of all subsets of
a set \(S\) is called the {\it powerset} of \(S\), and written as 
$2^{S}$. Cantor has shown (diagonal procedure)
that the powerset of \(S\) will always have greater cardinality
than the set \(S\), even when \(S\) is an infinite set.
An important consequence of this is that
we can without end construct bigger and bigger sets, 
\[2^{\aleph_0},~2^{2^{\aleph_0}},~2^{2^{2^{\aleph_0}}},~\cdots\]
and hence in set theory we cannot have a set which has the 
highest cardinality. A
disappointing consequence is that we cannot have a universal collection as
part of set theory and such a collection will always have to be
outside the set theory. 
One-to-one correspondence is the basis on 
which cardinality is decided, from which it follows that $\aleph_0$ 
can also be written as 
\[
\{1,2,4,8,\cdots\}=\{2^{0},2^{1},2^{2},2^{3},\cdots\}.
\]
As Halmos points out \cite{Hal:NST}, there is confusion 
in the literature regarding 
the notation $2^{\omega}$, it has been used to represent the 
above set and also the set $2^{\aleph_0}$, 
which \textit{in extenso}, can be written as 
\[
2^{\{{0,1,2,3,\cdots}\}}.
\]
To prevent this confusion, whenever necessary, we will write 
$2^{\aleph_0}$ as 
\[
\{<2^{0},2^{1},2^{2},2^{3},\cdots>\}
\]
to imply that $2^{\aleph_0}$ is a derived set from 
\[
\{2^{0},2^{1},2^{2},2^{3},\cdots\}.
\]

\section{INFINITESIMALS}

The study of the set of natural numbers gave us the notion 
of $\aleph_{0}$. The concept of a powerset makes it clear that there 
are higher cardinals above $\aleph_{0}$. Then the question arises, 
whether there is some other way of generating larger cardinals, 
other than taking powersets. Cantor has shown that this is 
possible, and gives us the sequence of transfinite sets of increasing 
cardinality as 
\[
\aleph_{0},\aleph_{1},\aleph_{2},\aleph_{3},\cdots ,
\]
with the understanding that there is no cardinal between $\aleph_{\alpha}$ 
and $\aleph_{\alpha+1}$. How exactly this sequence was generated, is 
an issue that we will take up later, but for the moment we will accept this 
sequence.

Because of the one-to-one correspondence between the right-sets and 
the left-sets, we will concentrate our attention on just the 
right-sets and right-sequences. Note, as an example, that the 
infinite sequence $.110****\cdots$ can be used
to represent the interval $(.75,.875)$, if we accept certain assumptions
about the representation:

\begin{quote}

The initial binary string, $.110=.75$, represents the initial point
of the interval. \\
The length of the binary string, $3$ in our case, decides the length
of the interval as $2^{-3}=.125$. \\
Every $*$ in the infinite $*$-string can be substituted by a
$0$ or $1$, to create $2^{\aleph_{\alpha}}$ points in the interval.

\end{quote}
Now, consider the right-sequence 
\[
.10101010 \cdots **** \cdots
\]
and the corresponding right-set
\[
\{1,3,5,7,\cdots\aleph_{0},\cdots\aleph_{\alpha}\}^{+}.
\]
If we can attach a meaning to this right-sequence, it can be
only this: it represents the number $.6666\cdots$ with an 
\textit{infinitesimal} attached to it, the cardinality of the
set of points inside the infinitesimal being
$2^{\aleph_{\alpha}}$.

\section{AXIOM OF FUSION}

The upshot of all our discussion so far is the following: The unit
interval $(0,1)$ is a set of infinitesimals with cardinality 
$\aleph_{0}$, with each infinitesimal representing a computable   
number. From the method we used in the construction of the 
infinitesimal, it will not be unreasonable, if we claim that the 
infinitesimal is an 
integral unit from which none of its $2^{\aleph_{\alpha}}$
elements can be removed. A set from which, the axiom of choice (AC) cannot 
remove an element, we will call a \textit{bonded set} and the elements in 
it \textit{figments}. If a set contains only bonded sets as its elements, 
then we will call it a \textit{class of bonded sets} 
or just \textit{bonded class}. We will use the term 
\textit{virtual cardinality} to refer to the cardinality of a bonded 
class. The set of all subsets of $\aleph_\alpha$ 
of cardinality $\aleph_\alpha$ we will symbolize as 
$\binom{\aleph_{\alpha}}{\aleph_{\alpha}}$.
Our saying so, will not, of course, make anything 
a fact, so we introduce an axiom called \textit{fusion}.

\begin{AF}
$(0,1)={\binom{\aleph_{\alpha}}{\aleph_{\alpha}}}
=R \times 2^{\aleph_{\alpha}}$, where $x \times 2^{\aleph_{\alpha}}$
is a  bonded set.
\end{AF}

\noindent The axiom of fusion says that $(0,1)$ is a class of bonded
sets, called infinitesimals. Further, the cardinality of each 
infinitesimal is $2^{\aleph_{\alpha}}$, and the \emph{virtual cardinality}
of $(0,1)$ is $\aleph_0$.

\begin{CBT}
${\binom{\aleph_{\alpha}}{\aleph_{\alpha}}}= 2^{\aleph_{\alpha}}$.
\end{CBT}

\begin{proof}
A direct consequence of the axiom of fusion is that
\[
2^{\aleph_{\alpha}} \le {\binom{\aleph_{\alpha}}{\aleph_{\alpha}}}.
\]
Since, ${\binom{\aleph_{\alpha}}{\aleph_{\alpha}}}$ is a subset of
$2^{\aleph_{\alpha}}$,
\[
{\binom{\aleph_{\alpha}}{\aleph_{\alpha}}} \le 2^{\aleph_{\alpha}},
\]
and the theorem follows.
\end{proof}

\section{EXPLOSIVE OPERATORS}

Halmos explains \cite{Hal:NST} the 
generation of {$\omega_1$}, the ordinal corresponding to
{$\aleph_1$} from {$\omega$} as given below. 

\begin{quote}
\hskip 1.5truecm  ...~In this way we get 
successively $\omega$, $\omega 2$, $\omega 3$, $\omega 4$,
$\cdots$. An application of the axiom of substitution yields something 
that follows them all in the same sense in which $\omega$ follows the 
natural numbers; that something is $\omega^2$. After that the whole
thing starts over again: $\omega^2+1$, $\omega^2+2$,
$\cdots$, $\omega^2+\omega$, $\omega^2+\omega+1$, $\omega^2+\omega+2$,
$\cdots$, $\omega^2+\omega 2$, $\omega^2+\omega 2+1$, $\cdots$,
$\omega^2$ $+\omega 3$, $\cdots$, $\omega^2+\omega 4$, $\cdots$, 
$\omega^2 2$, $\cdots$, $\omega^2 3$, $\cdots$, $\omega^3$,
$\cdots$, $\omega^4$, $\cdots$, $\omega^{\omega}$, $\cdots$,
$\omega^{(\omega^{\omega})}$, $\cdots$, 
$\omega^{(\omega^{(\omega^{\omega})})}$, 
$\cdots$ $\cdots$. The next one after all this is $\epsilon_0$; then come
$\epsilon_0+1$, $\epsilon_0+2$, $\cdots$, $\epsilon_0+\omega$, $\cdots$,
$\epsilon_0+\omega 2$, $\cdots$, $\epsilon_0+\omega^2$, $\cdots$,
$\epsilon_0+\omega^{\omega}$, $\cdots$, $\epsilon_02$, $\cdots$,
$\epsilon_0\omega$, $\cdots$, $\epsilon_0\omega^{\omega}$,
$\cdots$, $\epsilon_0^2$, $\cdots$ $\cdots$ $\cdots$.
\end{quote}
We want to write the essence of this quotation as terse as possible, for 
this purpose, we will first define \textit{explosive operators}.
For positive integers $m$ and $n$, we define an infinite sequence of
operators as follows.
\begin{align*}
m\otimes^{0} n
&=mn, \\
m\otimes^{k} 1
&=m, \\
m\otimes^k n
&=m\otimes^h [m\otimes^h [\cdots [m\otimes^h m]]], 
\end{align*}
where the number of $m$'s in the product is $n$ and
$h=k-1$. It is easy to see that 
\begin{align*}
m\otimes^{1} n
&=m^{n}, \\
m\otimes^{2} n
&=m^{m^{.^{.^{.^m}}}},
\end{align*}
where the number of $m$'s tilting forward is $n$.
We can continue to expand the operators in this fashion further, 
straining our currently available notations, but 
it is not relevant for us here. 
Note that these explosive operators are nothing but the 
well-known Ackermann functions.
We use these operators for symbolizing the transfinite
cardinals of Cantor.

\section{AXIOM OF MONOTONICITY}

Stripped of all verbal explanations, we can write the generation of 
$\omega_1$ as
\[
<0,1,2,\cdots \omega , \cdots \omega ^2,
\cdots \omega ^\omega ,
\cdots {}^\omega \omega ,
\cdots, \cdots, \cdots>
\] 
or in terms of the explosive operators as
\[
<0,1,2,\cdots \omega ,\cdots \omega \otimes^0 \omega,
\cdots \omega \otimes^1 \omega, \cdots \omega \otimes^2 \omega ,
\cdots ,\cdots ,\cdots>.
\]
Cantor has shown that the cardinality of $\omega \otimes^k \omega$ 
is $\aleph_0$ for all finite values of $k$, and hence it is not that we 
have a sequence here of increasing cardinality. Taking into account 
this fact, we assert that what the sequence means is that
\begin{align*}
\aleph_1
&= \{<0,1,2,\cdots \omega ,\cdots \omega \otimes^0 \omega,
\cdots \omega \otimes^1 \omega, \cdots \omega \otimes^2 \omega ,
\cdots >\}\\
&= \aleph_0 \otimes^{\{0,1,2,\cdots\}}\aleph_0\\
&= \aleph_0 \otimes^{\aleph_{0}} \aleph_0.
\end{align*}
Once this is accepted, a natural extension is that 
\[
\aleph_{\alpha+1}= \aleph_{\alpha} \otimes^{\aleph_{0}} \aleph_{\alpha}.
\]
An inspection of the explosive operators shows that $m\otimes^k n$ 
is a monotonically increasing function of $m,k,$ and $n$. Hence it 
will not be unreasonable to expect $m\otimes^k n$ to remain 
at least monotonically nondecreasing, when $m,k,$ and $n$ assume 
tranfinite cardinal values. Our saying all this, will not make it a 
fact, for that reason we state an axiom called axiom of 
\textit{monotonicity}.
Cantor always wanted his Continuum Hypothesis, $2^{\aleph_0}=\aleph_1$, to
be true in his set theory. We now introduce an axiom that 
accomplishes this, and
even more.

\begin{AM}
$\aleph_{\alpha+1}= \aleph_{\alpha} \otimes^{\aleph_0} \aleph_{\alpha}$,
and $2^{\aleph_\alpha}= 2 \otimes^{1} \aleph_\alpha$. Further,
if $m_1 \le m_2$, $k_1 \le k_2$, and
$n_1 \le n_2$, then $m_1 \otimes^{k_1} n_1 \le m_2 \otimes^{k_2} n_2$.
\end{AM}

\begin{CT}
$\aleph_{\alpha+1}= m\otimes^k\aleph_\alpha$ for finite $m>1,k>0$.
\end{CT}

\begin{proof}
A direct consequence of the axiom of monotonicity is that, for finite $m>1$
and $k>0$,
\[
2^{\aleph_\alpha}= 2~\otimes^1 \aleph_\alpha \le
m\otimes^k \aleph_\alpha
\le \aleph_\alpha \otimes^{\aleph_0} \aleph_\alpha =
\aleph_{\alpha+1}.
\]
When we combine this with Cantor's result
\[
\aleph_{\alpha+1} \le 2^{\aleph_\alpha},
\]
the theorem follows.
\end{proof}

\begin{GCH}
$\aleph_{\alpha+1}= 2^{\aleph_\alpha}$.
\end{GCH}

\begin{proof}
If we put $m=2$, $k=1$ in the Continuum Theorem, we get
\[
\aleph_{\alpha+1}= 2~\otimes^1 \aleph_\alpha= 2^{\aleph_\alpha},
\]
making GCH a theorem.
\end{proof}

\begin{UT}
All the three sequences
\[ \begin{array}{rrrrr}
\aleph_0, &\aleph_1,     &\aleph_2,     &\aleph_3,     &\ldots \\
\aleph_0, &2^{\aleph_0}, &2^{\aleph_1}, &2^{\aleph_2}, &\ldots \\
\aleph_0, &{\binom{\aleph_0}{\aleph_0}}, &{\binom{\aleph_1}{\aleph_1}},
&{\binom{\aleph_2}{\aleph_2}}, &\ldots 
\end{array} \]
represent the same series of cardinals.
\end{UT}

\begin{proof}
The axiom of monotonicity shows that the first two are the same,
and the axiom of fusion shows that the last two are same.
\end{proof}

\begin{AC}
Cartesian product of nonempty sets will always be nonempty, even 
if the product is of an infinite family of sets.
\end{AC}

\begin{proof}
GCH implies AC, and we have already proved GCH.
\end{proof}

\section{CONCLUSION}

A new concept that we have introduced in IST is that of a bonded 
set containing figments. It is somewhat like the concept of quarks 
in particle physics, where we know that they are there, but 
we cannot get one of them isolated. Figments can be very helpful 
in visualizing the space around us. If we call an infinitesimal with 
figments in it a \textit{white hole}, we can say that the finite part 
of our physical space is nothing but a tightly packed set of 
white holes. Since every irrational number has an infinitesimal 
attached with it, we can claim that every supernatural number has 
a \textit{black stretch} attached with it and the physical space 
beyond the finite part is a \textit{black whole} containing black 
stretches in it. 

IST visualizes an infinite recursive subset of positive integers as a 
number in the interval $(0,1)$, with a corresponding infinitesimal. 
This infinitesimal has in it all the transfinite sets containing the 
original recursive set. 

In measure theory, it has not been possible to date to construct a 
nonLebesgue measurable set without invoking the axiom of choice. 
IST does not allow figments to be picked up by the axiom of choice 
and for that reason, it would not be unreasonable to say that there 
are no nonLebesgue measurable sets in IST.

If we ignore figments, we can visualize the interval $(0,1)$ as a
set with virtual cardinality $\aleph_0$. As a consequence,
the Skolem paradox cannot be a serious problem in IST.

More than anything else, IST tells us to be realistic. It maintains 
that there are points we cannot touch, and that there are spaces we cannot 
reach. 

\bibliographystyle{amsplain}

\end{document}